\documentclass[preprint,12pt]{article}
\usepackage[latin9]{inputenc}
\usepackage{float}
\usepackage{amsmath}
\usepackage{amssymb}
\usepackage{graphicx}
\usepackage{esint}

\usepackage{amsmath,amstext,amsfonts,amssymb,amsbsy,array,graphicx}
\usepackage{stmaryrd}
\usepackage{float}
\usepackage{subfig}

\renewcommand{\div}{\text{div$~$}}
\newcommand{\Div}{\text{{\bf div}$~$}}

\newcommand{\bfv}{{\boldsymbol v}}

\newcommand{\bfsig}{{\boldsymbol\sigma}}

\newcommand{\bfe}{{\boldsymbol e}}
\newcommand{\bftheta}{{\boldsymbol \theta}}

\newcommand{\bfphi}{{\boldsymbol \phi}}

\newcommand{\bftwei}{{\boldsymbol \vartheta}}

\newcommand{\CP}{{\mathcal C}{\mathcal P}}

\newcommand{\bfpi}{{\boldsymbol \pi}}

\newcommand{\tr}{\text{tr }}
\newcommand{\bfI}{{\boldsymbol I}}
\newcommand{\bfu}{{\boldsymbol u}}

\newcommand{\triangulation}{\mathcal{T}}
\newcommand{\mcE}{\mathcal{E}}
\newcommand{\mcT}{\mathcal{T}}
\newcommand{\mcP}{\mathcal{P}}
\newcommand{\mcN}{\mathcal{N}}

\newcommand{\bfU}{{\boldsymbol U}}

\newcommand{\bfV}{{\boldsymbol \Gamma_h}}

\newcommand{\bfn}{{\boldsymbol n}}

\newcommand{\bfSig}{{\boldsymbol\Sigma}}
\newcommand{\bfeps}{{\boldsymbol\varepsilon}}
\newcommand{\bff}{{\boldsymbol f}}

\newcommand{\bfx}{{\boldsymbol x}}

\newcommand{\tn}{|\mspace{-1mu}|\mspace{-1mu}|}
\def\eref#1{{\rm (\ref{#1})}}

\begin{document}
\title{A posteriori error estimates for continuous/discon\-tinuous Galerkin
approximations of the Kirchhoff--Love buckling problem}
\author{Peter Hansbo$^a$, Mats G. Larson$^b$\\[4mm] \it\small$^a$Department of Mechanical Engineering, \\\it\small J\"onk\"oping University\\\it\small
SE-55111 J\"onk\"oping, Sweden \\\it\small $^b$Department of Mathematics and Mathematical Statistics\\\it\small Ume{\aa} University\\\it\small 
SE-901 87 Ume{\aa}, Sweden}
\date{}
\maketitle
\begin{abstract}
Second order buckling theory involves a one-way coupled coupled problem
where the stress tensor from a plane stress problem appears in an eigenvalue
problem for the fourth order Kirchhoff plate. In this paper we present an a
posteriori error estimate for the critical buckling load and mode corresponding
to the smallest eigenvalue and associated eigenvector. A particular feature of
the analysis is that we take the effect of approximate computation of the stress
tensor and also provide an error indicator for the plane stress problem. The Kirchhoff
plate is discretized using a continuous/discontinuous finite element method which uses
standard continuous piecewise polynomial finite element spaces which can also be
used to solve the plane stress problem.
\end{abstract}

\section{Introduction}

Buckling of thin plates can be modeled by an eigenvalue problem involving the stress 
tensor of the plane stress problem corresponding to a given load situation tangential 
to the plate. The smallest eigenvalues corresponds to the critical parameter multiplying 
the given plane stress load that results in buckling.

Thin plates are modeled by fourth order differential equations according to the
Kirchhoff-Love theory and require special attention when discretized using the finite element
method. In this paper we use the continuous/discontinuous Galerkin (c/dG) method proposed
by Engel et al. \cite{EngGar02} which is based on standard continuous piecewise polynomial
spaces of order greater or equal to two inserted into a discontinuous Galerkin formulation,
see Hansbo and Larson \cite{HanLar02}, of the fourth order plate equation. We refer also to Wells and Dung \cite{WeDu07} for
a method closely related to the one presented here, and to Noels and Radovitzky \cite{Noels08} for an extension of the c/dG idea to Kirchhoff--Love shells.
 
The c/dG formulation 
has the advantage that it uses standard finite element spaces, is easy to implement, and
extends naturally to higher order polynomials. Another important advantage in this particular
problem is that we may solve the plane stress problem using the same finite element spaces. Note
that this would not be the case if we, for instance, used nonconforming Morley elements for the 
plate problem since these element can not be used for the plane stress problem.

In this paper we derive a posteriori error estimates for the critical buckling load and mode
corresponding to the first eigenpair, and use these estimates to obtain mesh refinement strategies for error reduction. The error estimates are derived using duality techniques and
are based on Larson \cite{Lar00} where a posteriori error estimates for the Poisson 
equation were presented. A particular feature of the estimates presented herein is that we also take the effect
of discretization of the plane stress problem into account. The error analysis of the buckling 
problem results in a specific goal functional which should be controlled in the plane stress solver. 
Here we follow the general approach to error estimation for one-way coupled problems developed by 
Larson and Bengzon \cite{LarBen07}, and adapted to linear second order plate theory in \cite{HaHeLa}. We also mention the work \cite{HeuRanXX} by Heuveline and Rannacher 
where a posteriori error estimates for a nonsymmetric eigenvalue problem related to the linearized 
stability of the Navier-Stokes equations is presented. These estimates also involve the effect of 
the accuracy in the computed flow field on the eigenvalue problem and are thus related to our approach.

This paper is organized as follows: in Section 2 we pre\-sent the Kirchhoff-Love
buckling problem and the continu\-ous/dis\-conti\-nuous Galerkin method, in Section 3 we
derive the {\em a posteriori}\/ error estimates, in Section 4 we present some
numerical results, and in Section 5 we present some conclusions.

\section{The Buckling Problem and Finite Element Method}

\subsection{The Kirchhoff-Love Buckling Eigenvalue Problem}

The clamped Kirchhoff--Love buckling problem takes the form: find the plate 
displacements $u_P$ (orthogonal to the plate) such that
\begin{align}
\div \Div \bfsig_P( \nabla u_P ) - \div t (\bfsig_M \nabla u_P) = f_P &  \text{ in $\Omega$}
\\
u_P = 0  &\text{ on $\partial \Omega$}
\\
\bfn \cdot \nabla u_P = 0   &\text{ on $\partial \Omega$}
\end{align}
where $t$ denotes the thickness of the plate and 
\begin{equation}\label{sigm}
\bfsig_M = 2 \mu \bfeps(\bfu_M) + \lambda \text{\tr}\bfeps(\bfu_M)\bfI
\end{equation}
where $\mu$ and $\lambda$ are the Lam\'e parameters, tr denotes the trace operator and $\bfI$\/ is the identity matrix, is determined 
by the membrane equation: find the membrane displacements $\bfu_M$ (tangential to the plate) such that
\begin{alignat}{3}
- \Div \bfsig_M (\bfu_M) &= \bff_M & \quad &\text{in $\Omega$}
\\
\bfu_M &= 0 &  &\text{on $\partial \Omega$}
\end{alignat}
Here
\begin{align}
\bfsig_P(\bfeps) = \frac{Et^3}{12(1-\nu^2)} \left( (1-\nu) \bfeps + {\nu}
\text{tr} ( \bfeps ) \, \bfI \right)
\end{align}
is the plate stress tensor, $\bfeps(\bfv) = (\nabla \bfv + (\nabla \bfv)^T)/2$ is the 
strain tensor, $E$ is the Young's modulus, and $\nu$ is the Poisson ratio, in terms of which 
$\lambda=E/(1+\nu)$ and $\mu=E\nu/(1-\nu^2)$.

Scaling the membrane forcing by a parameter $\lambda_P$, i.e., replacing the load by
$\lambda_P \bff_M$ we note that by linearity $\bfsig_M$ is replaced by $\lambda_P \bfsig_M$.
The critical buckling loads are then determined by the eigenvalue problem: find $u_P$ and $\lambda_P$
such that
\begin{align}
\div \Div \bfsig_P( \nabla u_P ) - \div t (\lambda_P \bfsig_M \nabla u_P) =  0  &  \text{ in $\Omega$}
\\
u_P = 0 & \text{ on $\partial \Omega$}
\\
\bfn \cdot \nabla u_P = 0 & \text{ on $\partial \Omega$}
\end{align}

The corresponding variational formulation reads: find the plate displacement $u_P \in
H^2_0(\Omega)$ and eigenvalue $\lambda_P \in {\bf R}$ such that
\begin{equation}\label{eq:varform}
a_P(\nabla u_P, \nabla v ) + \lambda_P (\bfsig_M \nabla u_P, \nabla v ) = 0 \quad \forall v \in H^2_0(\Omega)
\end{equation}
where $ \bfsig_M\in [L^2(\Omega)]^{2 \times 2} $ defined by (\ref{sigm}),
with $\bfu_M \in [H^1_0(\Omega)]^2$ the solution of
\begin{equation}\label{eq:varform}
a_M(\bfu_M, \bfv ) = (\bff_M, \bfv) \quad \forall \bfv \in [H^1_0(\Omega)]^2
\end{equation}
Here the bilinear forms $a_P(\cdot,\cdot)$ and $a_M(\cdot,\cdot)$ are defined by
\begin{align}
a_P(\bftheta,\bftwei) & = (\bfsig_P(\bftheta), 
\bfeps(\bftwei))  \\ a_M(\bftheta,\bftwei) & = (\bfsig_M(\bftheta), \bfeps(\bftwei))
\end{align}
where $(\cdot,\cdot)$ is the $L^2(\Omega)$ inner product.

\subsection{The Mesh and Finite Element Spaces}

We consider a subdivision $ \triangulation =\{ T\}$ of $\Omega$ into
a geometrically conforming finite element mesh. We assume that the
elements are shape regular, i.e., the quotient of the diameter of the
smallest circumscribed sphere and the largest inscribed sphere is uniformly
bounded. We denote by $h_T$ the diameter of element $T$ and by
$h = \max_{T \in \mcT} h_T$ the global mesh size parameter. We shall use
continuous, piecewise polynomial, approximations of the transverse
displacement:
\begin{equation}
\CP_{k} =\{ v \in C^0(\Omega):~\text{$v \vert_T \in \mcP_{k}(T)$
$\forall T \in \triangulation$}\}
\end{equation}
where $\mcP_k(T)$ is the space of polynomials of order $k \geq 2$
defined on $T$. Furthermore, we let $\CP_{k,0} = \CP_{k} \cap H^1_0$.

We introduce the Scott-Zhang interpolation operator $\pi: H^1_0(\Omega)
\rightarrow \CP_{k,0}$ and recall the following elementwise interpolation
error estimate
\begin{equation}\label{interpol}
| u - \pi u |_{T,m} \leq C h_T^{s-m}| u |_{\mcN(T),s}
\end{equation}
where $0\leq m \leq s \leq k+1$ and $\mcN(T)$ is the union of all elements
which are neighbors to element $T$.

To define our method we introduce the set of edges in the mesh,
$\mcE =\{ E \}$, and we split $\mcE$ into two disjoint subsets
\begin{equation}
  \mcE = \mcE_I \cup \mcE_B
\end{equation}
where $\mcE_I$ is the set of edges in the interior of $\Omega$ and
$\mcE_B$ is the set of edges on the boundary. Further, with each
edge we associate a fixed unit normal $\bfn$  such that for edges
on the boundary $\bfn$ is the exterior unit normal. We denote the
jump of a function $\bfv \in \bfV$ at an edge $E$ by $ \left[\bfv
\right] = \bfv^+-\bfv^-$ for $E \in \mcE_I$ and $
\left[\bfv\right] = \bfv^+$ for $E \in \mcE_B$, and the average
$\langle  \bfv\rangle = (\bfv^+ + \bfv^-)/2$ for $E \in \mcE_I$
and $\langle \bfv\rangle = \bfv^+$ for $E \in \mcE_B$, where
$\bfv^{\pm} = \lim_{\epsilon\downarrow 0} \bfv(\bfx\mp
\epsilon\,\bfn)$ with $\bfx\in E$.

\subsection{The Continuous/Discontinuous Galerkin Method}

We shall solve the membrane equation using standard continuous Galerkin and
the plate problem with the continuous/discontinuous Galerkin method. The method takes the form: find $U_P \in \CP_{k_P,0}$ and $\Lambda_P \in {\bf R}$
such that
\begin{equation}\label{eq:cdg}
A_{P}(\nabla U_P, \nabla v ) + \Lambda_P t (\bfSig_M \nabla U_P, \nabla v) = 0 \quad\forall v \in \CP_{k_P,0} 
\end{equation}
where  $ \bfSig_M = 2 \mu \bfeps(\bfU_M) + \lambda \text{\tr}\bfeps(\bfU_M)  \bfI$
with $\bfU_M \in  [\CP_{k_M,0}]^2$ determined by
\begin{equation}\label{eq:gom}
a_M(\bfU_M, \bfv ) = (\bff_M, \bfv) \quad \forall \bfv \in  [\CP_{k_M,0}]^2
\end{equation}
The bilinear form $A_P(\cdot,\cdot)$ is defined by
\begin{align} \nonumber
A_P(\bftheta,\bftwei) = {}& \sum_{T \in \mcT} ( \bfsig(\bftheta),
\bfeps( \bftwei))_T
\\ \nonumber
&  - \sum_{E \in \mcE_I \cup \mcE_B }
 ( \langle\bfn \cdot \bfsig(\bftheta) \rangle,  [ \bftwei
  ])_E \\ \nonumber
& -  \sum_{E \in \mcE_I \cup \mcE_B }( [\bftheta ], \langle \bfn \cdot \bfsig(\bftwei)\rangle)_E
\\ 
&  + (2 \mu +2\lambda)\, \gamma \sum_{E \in \mcE_I \cup
\mcE_B}
         h_E^{-1} ([ \bftheta ], \left[ \bftwei\right])_E\label{nitsche_form}
\end{align}
for all $\bftheta,\bftwei \in \bigoplus_{T \in \mcT} [H^1(T)]^2$.
Here $\gamma$ is a positive parameter and $h_E$ is defined by
\begin{equation}
h_E = \left( |T^+| + |T^-| \right) / ( 2 \, |E| ) \quad \text{for $E
= \partial T^+ \cap
\partial T^-$}
\end{equation}
with $|T|$ the area of $T$, on each edge $E$. See \cite{HanLar09} for
details on the value of $\gamma$.

\section{A Posteriori Error Estimates}

\subsection{Preliminaries}

We first define a projector onto $E(\lambda_P)$ that is associated with the natural scalar 
products involved in the variational statement. We define $\mcP_{\lambda_P}: H^1_0(\Omega) 
\rightarrow E(\lambda_P)$ as follows
\begin{equation}
A_P(\nabla \mcP_{\lambda_P} v, \nabla w ) = A_P( \nabla v, \nabla w ) \quad \forall w \in E(\lambda_P)
\end{equation}
Note that since $w$ are eigenfunctions associated with $\lambda_P$ the projection 
also satisfies the following equation 
\begin{equation}
(\sigma_M \nabla \mcP_{\lambda_P} v , \nabla w) =  (\sigma_M \nabla v, \nabla w) \quad \forall w \in E(\lambda_P)
\end{equation}
We introduce the norm 
\begin{equation}
\tn v \tn^2 = A_P( \nabla v, \nabla v ), \quad \forall v \in H^2_0(\Omega) \cup \CP_{k_P,0}
\end{equation}
and normalize computed eigenvectors $U_P$ as follows 
\begin{equation}
\tn \nabla U_P \tn^2  = A_P(\nabla U_P, \nabla U_P) 
 = |\Lambda_P (\Sigma_M \nabla U_P, \nabla U_P )| = 1\label{scaling}
\end{equation}



\subsection{Error Representation Formulas}

\paragraph{The dual problem.}

To derive error representation formulas we introduce the following dual problem: find $\phi_P$ such that
\begin{align}
\div \Div \bfsig_P( \nabla \phi_P ) - \div t \lambda_P (\bfsig_M \nabla \phi_P) &= \psi_P  \; \text{in $\Omega$}
\\
\phi_P &= 0 \;\text{on $\partial \Omega$}
\\
\bfn \cdot \nabla \phi_P &= 0 \;\text{on $\partial \Omega$}
\end{align}
Different choices of the righthand side will lead to estimates for the errors in eigenvalues and eigenvectors.
The righthand sides $\psi_P$ will be chosen in such a way that the solution to the dual problem is well defined. We return to these issues below.

Multiplying with the error $e_P = u_P - U_P$ and integrating by parts we obtain
\begin{align}
(e_P,\psi_P) = {}& (e_P, \div \Div \bfsig_P( \nabla \phi_P ))\nonumber \\\nonumber
 & - (e_P,\div t \lambda_P (\bfsig_M \nabla \phi_P))\\\nonumber
= {}& \sum_{T \in \mcT} ( \bfsig_P( \nabla e_P ), \bfeps(\nabla \phi_P) )_T \\\nonumber
& - \sum_{E \in \mcE_I \cup \mcE_B} ([\nabla e_P ], \bfn \cdot\bfsig_P(\nabla \phi_P ) )_E
\\ \nonumber
&  - \lambda_P t (\bfsig_M \nabla U_P, \nabla \phi_P )
\\ 
& - A_P( \nabla U_P, \nabla \phi_P  ) - \lambda_P t (\bfsig_M \nabla U_P, \nabla \phi_P )\label{eq:errorrepa}
\\  \nonumber
={}& - A_P( \nabla U_P, \nabla (\phi_P - \pi_P \phi_P)  ) \\\nonumber
& - \Lambda_P t (\bfSig_M \nabla U_P, \nabla (\phi_P - \pi_P \phi_P) )
\\ \nonumber
& + (\Lambda_P - \lambda_P)  t (\bfsig_M \nabla U_P, \nabla \phi_P  )
\\ 
& -  \Lambda_P t ( (\bfsig_M - \bfSig_M) \nabla U_P, \nabla \phi_P )\label{eq:errorrepb}
\end{align}
where in \eref{eq:errorrepa} we used the fact that $[\nabla \phi_P ] = 0$ and then in \eref{eq:errorrepb}
we rearranged the terms using the identity 
$\lambda_P \bfsig_M = \Lambda_P \bfSig_M - (\Lambda_P - \lambda_P) \bfsig_M
+ \Lambda_P (\bfsig_M - \bfSig_M)$
and finally used Galerkin orthogonality \eref{eq:cdg} to subtract $\pi_P \phi_P$.

\paragraph{Representation of the Error in the Eigenvalue}

Setting $\psi_P = 0$ and denoting the solution to the dual problem by $\phi_{P,\lambda_P}$ we get
\begin{align}
& (\Lambda_P - \lambda_P)  t (\bfsig_M \nabla U_{P}, \nabla \phi_{P,\lambda_P}  )  \nonumber
\\ \nonumber
& \quad =A_P( \nabla U_P, \nabla (\phi_{P,\lambda_P} - \pi_P \phi_{P,\lambda_P})  )\\ \nonumber
&     \qquad          + \Lambda_P t (\bfSig_M \nabla U_P, \nabla (\phi_{P,\lambda_P} - \pi_P \phi_{P,\lambda_P}) )
\\  \label{eq:errorrepeig_a}
&     \qquad  +  \Lambda_P t ( (\bfSig_M - \bfsig_M) \nabla U_P, \nabla \phi_{P,\lambda_P} )&
\end{align}
In this case the solution to the dual problem is an arbitrary eigenfunction associated with $\lambda_P$, 
i.e. $\phi_{P,\lambda_P} \in E(\lambda_P)$. Choosing $\phi_{P,\lambda_P} = \mcP_{\lambda_P} U_P/ \tn \mcP_{\lambda_P} U_P \tn $ we obtain the following estimate
\begin{align} \nonumber
& | (\Lambda_P - \lambda_P)  t (\bfsig_M \nabla U_P, \nabla \phi_P  ) |  \\ \nonumber
& \quad = | (\Lambda_P - \lambda_P) \lambda_P^{-1}  A_P( \nabla U_P, \nabla \phi_{P,\lambda_P}  ) |
\\ \nonumber
&\quad=| (\Lambda_P - \lambda_P) \lambda_P^{-1} | A_P( \nabla \mcP_{\lambda_P} U_P, \nabla \phi_{P,\lambda_P}  ) |
\\
&\quad= | (\Lambda_P - \lambda_P) \lambda_P^{-1} | \, \tn \nabla \mcP_{\lambda_P} U_P \tn\label{esteigbb}
\end{align}
We now assume that the computed eigenvalue $\Lambda_P$ approximates
the exact eigenvalue $\lambda_P$ and that there are constants $0 \leq \delta < 1$
and $h_0$ such that
\begin{equation}\label{assumptiona}
\tn \nabla (I - \mcP_{\lambda_P} ) U_P \tn \leq \delta
\end{equation}
for all meshes with $ \max_{T \in \mcT} h_T \leq h_0 $. We remark that the validity of this 
assumption follows from standard a priori convergence theory. Using \eref{esteigbb}, 
\eref{assumptiona}, and the scaling \eref{scaling} together with Pythagoras identity we obtain 
\begin{equation}
| (\Lambda_P - \lambda_P)  t (\bfsig_M \nabla U_P, \nabla \phi_P  ) |
  \label{eq:errorrepeig_b} \geq | (\Lambda_P - \lambda_P) \lambda_P^{-1} | ( 1 - \delta^2 )^{1/2}
\end{equation}
Finally, combining \eref{eq:errorrepeig_a}, \eref{eq:errorrepeig_b}, and using 
the triangle inequality we arrive at 
\begin{align}
&(1- \delta^2)^{1/2} | (\Lambda_P - \lambda_P)\lambda_P^{-1} |  \nonumber 
\\ \nonumber
& \quad \leq | A_P( \nabla U_P, \nabla (\phi_{P,\lambda_P} - \pi_P \phi_{P,\lambda_P})  )
    \\ \nonumber &\qquad             + \Lambda_P t (\bfSig_M \nabla U_P, \nabla (\phi_{P,\lambda_P} - \pi_P \phi_{P,\lambda_P}) ) |
\\
& \qquad         +  | \Lambda_P t ( (\bfsig_M - \bfSig_M) \nabla U_P, \nabla \phi_{P,\lambda_P} ) |
\end{align}

\paragraph{Representation of the Error in the Eigenvector.}

Following Larson \cite{Lar00} we define the error in an eigenvector to be the component
orthogonal to the exact eigenspace which it approximates an element in. Note that this 
definition has the advantage that it covers also multiple eigenvectors. More precisely we 
will estimate the error in the $H^m(\Omega)$ seminorm for $m=0,1$. We then define the error 
$e_m$ as 
\begin{equation}
e_m = (I - P_m ) U_P 
\end{equation}
where $P_m$ is the orthogonal projection $H^m(\Omega) \rightarrow E(\lambda_P)$ defined by 
$(v - P_0 v, w ) = 0$ and $(\nabla (v - P_1 v, \nabla w )$ for all $w\in E(\lambda_P)$ and 
$m=0,1,$ respectively. To represent the  semi norm $| e_P |_m$ we let 
$\psi_P = \psi_{P,m} = (-\Delta)^m e_{P,m} / | e_{P,m} |_m $ with $m=0,1$ and we denote the 
corresponding solution to the dual problem by $\phi_{P,m}, m=0,1$. We then get
\begin{align}\nonumber
| e_{P,m} |_m ={}&- A_P( \nabla U_P, \nabla (\phi_{P,m} - \pi_P \phi_{P,m})  ) \\ \nonumber
& - \Lambda_P t (\bfSig_M \nabla U_P, \nabla (\phi_{P,u} - \pi_P \phi_{P,m}) )
\\ \nonumber
& + (\Lambda_P - \lambda_P)  t (\bfsig_M \nabla U_P, \nabla \phi_{P,m}  )
\\ \label{eq:eigvector_a}
& +  \lambda_P t ( (\bfSig_M - \bfsig_M) \nabla U_P, \nabla \phi_{P,m} )
\end{align}
In this case we require the solution $\phi_{P,m}$ to be orthogonal to $E(\lambda_P)$
to achieve uniqueness. 

Next we estimate the second term on the right hand side as follows
\begin{align} \nonumber
\lambda_P t (\bfsig_M \nabla U_P, \nabla \phi_{P,m}  ) &= A_P (\nabla U_P, \nabla \phi_{P,m})
\\ \nonumber
&= A_P (\nabla e_{P,m}, \nabla \phi_{P,m})
\\ \nonumber
&\leq | e_{P,m}|_m | \phi_{P,m} |_{4-m}
\\ \nonumber
& \leq C_m | e_{P,m} |_m | \psi_{P,m} |_{-m}
\\
& \leq  C_m | e_{P,m} |_m 
\end{align}
where we used the stability estimate $$| \phi_{P,m} |_{4-m} \leq C |\psi_{P,m} |_{-m}$$ 
and at last the identity $| \psi_{P,m} |_{-m} = 1$ which follows from the definition of 
$\psi_{P,m}$. Thus we have
\begin{align}
& | (\Lambda_P - \lambda_P)  t (\bfsig_M \nabla U_P, \nabla \phi_{P,m}  ) | \nonumber\\ \label{estaa}
&\quad\leq | (\Lambda_P - \lambda_P) \lambda_P^{-1} | C_m | e_{P,m} |_m 
\end{align}
Now again assuming that the computed eigenvalue $\Lambda_P$ approximates
the exact eigenvalue $\lambda_P$ and that there are constants $0 \leq \delta < 1$
and $h_0$ such that
\begin{equation}\label{assumptionb}
| (\Lambda_P - \lambda_P) \lambda_P^{-1} | C_m \leq \delta
\end{equation}
for all meshes with $ \max_{T \in \mcT} h_T \leq h_0 $. We note again that the validity 
of this assumption follows from standard a priori convergence theory.
Combining (\ref{eq:eigvector_a}), (\ref{estaa}), and (\ref{assumptionb}) and using the 
triangle inequality we obtain the estimate
\begin{align} \nonumber
 (1 - \delta ) | e_P |_m \nonumber \leq {}& | A_P( \nabla U_P, \nabla (\phi_{P,m} - \pi_P \phi_{P,m})  ) \\ \nonumber
& + \Lambda_P t (\bfSig_M \nabla U_P, \nabla (\phi_{P,m} - \pi_P \phi_{P,m}) ) |
\\
& +  | \lambda_P t ( (\bfSig_M - \bfsig_M) \nabla U_P, \nabla \phi_{P,m} ) |
\end{align}
{\bf Remark.} The constant $C_m$ is of the form
\begin{equation}
C_m = \frac{c_m}{\text{gap}(\lambda_P)}
\end{equation}
where $\text{gap}(\lambda_P)$ is the distance between $\lambda_P$ and the closest eigenvalue. 
Thus assumption \eref{assumptionb} guarantees satisfactory resolution of the spectrum in the 
vicinity of $\lambda_P$.

\paragraph{Representation of the Modeling Error.}

Introducing the dual problem: find $\bfphi_{M,X} \in [H^1_0(\Omega)]^2$ such that
\begin{equation}
a_M(\bfv,\bfphi_{M,X}) = \lambda t (\bfsig_M(\bfv)\nabla U_P,\nabla \phi_{P,X}) ,
\end{equation}
for all $\bfv \in [H^1_0(\Omega)]^2, X \in \{0,1,\lambda_P\}$,
we get, by setting $\bfv =\bfe_M:=\bfu_M-\bfU_M$ and using Galerkin orthogonality \eref{eq:gom} for the membrane equation,
the following error representation formula
\begin{align}\nonumber
& ( (\bfSig_M - \bfsig_M) \nabla U_P, \nabla \phi_{P,X} ) \nonumber \\
& \quad = (\bfsig_M(\bfe_M)\nabla U_P,\nabla \phi_{P,X})\nonumber
\\
& \quad=a_M(\bfe_M,\bfphi_{M,X})\nonumber
\\\nonumber
& \quad= a_M(\bfe_M,\bfphi_{M,X} - \bfpi_M \bfphi_{M,X})
\\
& \quad = (\bff_M, \bfphi_{M,X} - \bfpi_M \bfphi_{M,X}) - a_M(\bfU_M,\bfphi_{M,X} - \bfpi_M \bfphi_{M,X})
\end{align}

\subsection{Abstract A Posteriori Error Estimates}

Combining the estimates above we obtain the following abstract error estimates. For 
the error in the eigenvalue
\begin{align}
& (1- \delta^2)^{1/2} |( \Lambda_P - \lambda_P ) \lambda_P^{-1} | \nonumber
\\\nonumber
&\quad \leq \vert A_P( \nabla U_P, \nabla (\phi_{P,\lambda_P} - \pi_P \phi_{P,\lambda_P})  )\\ \nonumber
     &\qquad            + \Lambda_P t (\bfsig_M \nabla U_P, \nabla (\phi_{P,\lambda_P} - \pi_P \phi_{P,\lambda_P}) ) \vert
\\ \nonumber
& \qquad  +   | (\bff_M, \bfphi_{M,\lambda_P} - \bfpi_M \bfphi_{M,\lambda_P}) \\ 
&\qquad-  a_M(\bfU_M,\bfphi_{M,\lambda_P} - \bfpi_M \bfphi_{M,\lambda_P})  |
\end{align}
and for the error in the eigenvector
\begin{align}\nonumber
(1 - \delta ) | e_P |_m \leq {}& | A_P( \nabla U_P, \nabla (\phi_{P,m} - \pi_P \phi_{P,m})  ) \\ \nonumber
& + \Lambda_P t (\bfSig_M \nabla U_P, \nabla (\phi_{P,m} - \pi_P \phi_{P,m}) ) |
\\ \nonumber
& + | (\bff_M, \bfphi_{M,m} - \bfpi_M \bfphi_{M,m}) \\
& -  a_M(\bfU_M,\bfphi_{M,m} - \bfpi_M \bfphi_{M,m}) |
\end{align}
for $m = 0,1$.

\subsection{Error Estimates Using the Dual Weighted Residual Approach}

Using standard procedures, involving integration by parts, the Cauchy-Schwartz inequality,
a trace inequality, and the interpolation error estimate (\ref{interpol}), we obtain the 
following estimate
\begin{align}\nonumber
& |A_P( \nabla U_P, \nabla (\phi_{P,X} - \pi_P \phi_{P,X})  )  
+  \\ \nonumber
&\qquad \Lambda_P t (\bfSig_M \nabla U_{P}, \nabla (\phi_{P,X} - \pi_P \phi_{P,X}) ) | \\
& \quad\leq \sum_{T \in \mcT_h} R_{P,T} W_{P,X,T}
\end{align}
where the plate element residual $R_{P,T}$ and weight $W_{P,X,T}$ are defined by
\begin{align}\nonumber
R^2_{P,T}= {}&\|f_P - \div \Div \bfsig_P(\nabla U_P)\|_T^2 \\\nonumber
& + h_T^{-1}
\| [\bfn \cdot \Div \bfsig_P ( \nabla U_P )] \|^2_{\partial T}
\\ 
& + h_T^{-3} \| [\bfn \cdot \bfsig_P ( \nabla U_P )] \|^2_{\partial
T}\nonumber \\ & + \gamma^2 h_E^{-5} \| [\nabla U_P ] \|^2_{\partial T}
\end{align}
\begin{equation}
W_{P,X,T} = h_T^{\alpha_{P,X}} | \phi_{P,X} |_{\mcN(K),\alpha_P}, \quad 0 \leq \alpha_{P,X} \leq k_P + 1
\end{equation}
Here the regularity parameter $\alpha_{P,X}$ reflects the regularity properties 
of the solutions to the dual problems. For the membrane problem we have the 
corresponding estimate
\begin{equation}
|a_M(\bfe_M,\bfphi_{M,X} - \bfpi_M \bfphi_{M,X})|
\leq \sum_{T \in \mcT_h} R_{M,T} W_{M,X,T}
\end{equation}
where the residual and weight are defined by
\begin{align}\nonumber
R^2_{M,T} ={}& \|\bff_M + \Div \bfsig_M(\bfU_M)\|_T^2 \\
& + h_T^{-1}
\| [\bfn \cdot \Div \bfsig_M ( \bfU_M )] \|^2_{\partial T}
\end{align}
and
\begin{equation}
W_{M,X,T} =  h_T^{\alpha_M} | \phi_{M,X} |_{\mcN(K),\alpha_{M,X}} 
\end{equation}
for $0 \leq \alpha_{M,X} \leq k_M + 1, X \in \{\lambda_P,0,1\}$.
Collecting these estimates and the abstract a posteriori error estimates 
we finally arrive at the following dual weighted residual a posteriori 
error estimates
\begin{align}\nonumber
(1- \delta^2)^{1/2} |(\Lambda_P - \lambda_P)\lambda_P^{-1} | 
\leq {}&\sum_{T \in \mcT_h} R_{P,T} W_{P,\lambda, T} \\
& +  \sum_{T \in \mcT_h} R_{M,T} W_{M,\lambda,T}
\end{align}
and for the error in the eigenvector
\begin{align}
(1 - \delta ) | e_P |_{m} & \leq \sum_{T \in \mcT_h} R_{P,T} W_{P,m,T} +  \sum_{T \in \mcT_h} R_{M,T} W_{M,m,T}
\end{align}
$m=0,1$. Considering the expected optimal regularity of the dual problems we may expect
\begin{equation}
\alpha_{P,\lambda} = k_P + 1, \quad \alpha_{P,m} = 4-m, \quad \alpha_{M,m} = \alpha_{M,\lambda} = 2
\end{equation}

\subsection{Residual Based Estimates}

Using stability estimates for the solutions to the dual problems we obtain the residual based estimates
\begin{align}\nonumber
& (1- \delta^2) |(\Lambda_P - \lambda_P)\lambda_P^{-1} |^2 \\
&\quad \leq  C\left( \sum_{T \in \mcT_h} h_T^{2\alpha_{P,\lambda_P}} R^2_{P,T} 
 + \sum_{T \in \mcT_h} h_T^{2\alpha_{M,\lambda}} R^2_{M,T}\right)
\end{align}
and for the error in the eigenvector
\begin{equation}
(1 - \delta )^2 | e_P |^2_{m} \leq C \left( \sum_{T \in \mcT_h} h_T^{2\alpha_{P,m}} R^2_{P,T} 
+ \sum_{T \in \mcT_h} h_T^{2\alpha_{M,m}} R^2_{M,T}\right)
\end{equation}
for $m=0,1$.

\section{Numerical examples}

\subsection{Known stress tensor}

We consider the L--shaped domain $$\Omega := (0,1)\times(0,1)\setminus (1/2,1)\times (0,1/2)$$
The plate is simply supported on all boundaries ($u=0$), and the in-plane stress tensor is
chosen as the unit tensor. 
Thus, we have no error contribution from the membrane problem. We set $E=1$, $\nu=1/4$, and $t=1$.
We use the adaptive algorithm for the computation of the
lowest three eigenvalues. The singularity in the inward-pointing corner is excited for the
first two but not for the third, which is also clearly visible in the adaptation of
the meshes shown in Figures 1--3. In Figure 4, we give the corresponding eigensolution, and in Figure 5 we give the corresponding effectivity indices (approximate error in eigenvalue divided by exact error). The third eigenvalue can be computed analytically, the first two have been estimated by an approximate solution on a dense mesh. The effectivity indices have been computed on a sequence of meshes obtained using a fixed ratio refinement technique where 
the elements with the highest 25\% element error indicators have been refined in each step.
The unknown constant in the error representation formula has been set so that the effectivity index is of medium size; the same constant has been used for all three eigenvalues.

\subsection{Computed stress tensor}

For our second example, we use the same domain, material data, and boundary conditions for the plate. For the 
elasticity computations, we use a body force $\bff = (r,-9\,r/10)$, where $r$ denotes the distance from the inward pointing corner.
The boundary conditions were: clamped conditions at $x=1/2$, $y\leq 1/2$, at $y=0$, at $x=1$, and at $y=1/2$, $x\geq 0$. The remaining boundaries were traction free. 

In Figure 5 we give the adapted mesh using the full estimate, and, for comparison, we also give, in Figures 6--7, the corresponding meshes 
when only partial estimates, plate residual and stress residual, respectively, are used. In Figure 9 we show the lowest buckling mode for which the estimate is aiming. Finally, we show, in Figure 10, how the different residuals behave asymptotically as estimates of the eigenvalue error. Clearly, in order to obtain an effectivity index
that does not increase or decrease, we need the full residual, though we concede that the balance between these two residuals may
be difficult to ascertain. We have here willfully chosen the balance in order to obtain a reasonably constant effectivity index for the full residual.

\section{Conclusions}

We have formulated a continuous/discontinuous Ga\-lerkin method for the thin plate buckling problem. 
The method has the advantage that we can solve both the membrane and plate problem 
with the same standard finite element spaces of continuous piecewise polynomials defined on triangles 
or quadrilaterals. Furthermore, we proved a posteriori error estimates for both the error in the eigenvalue 
(critical buckling load) and the eigenvectors (buckling modes) with the special feature that also 
the effect of approximate solution of the membrane problem is taken into account. Based on the estimates 
we constructed an adaptive algorithm for adaptive mesh refinement.

\newpage 

\begin{figure}
 \begin{minipage}[c]{0.42\textwidth}
\includegraphics[width=4in]{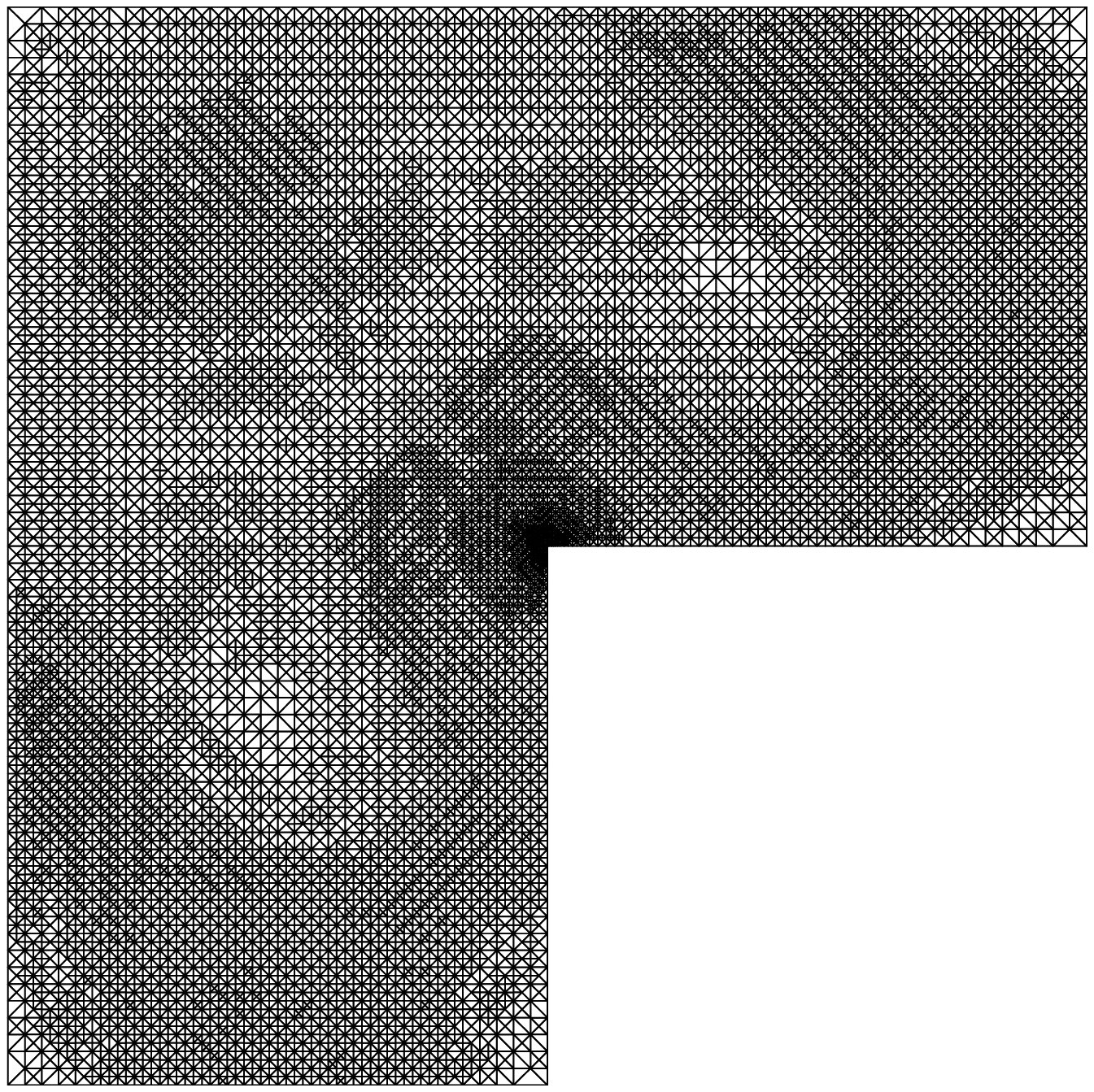}
 \end{minipage}
\caption{Adapted mesh  for the first eigenvalue}
\end{figure}
\begin{figure}
 \begin{minipage}[c]{0.42\textwidth}
\includegraphics[width=4in]{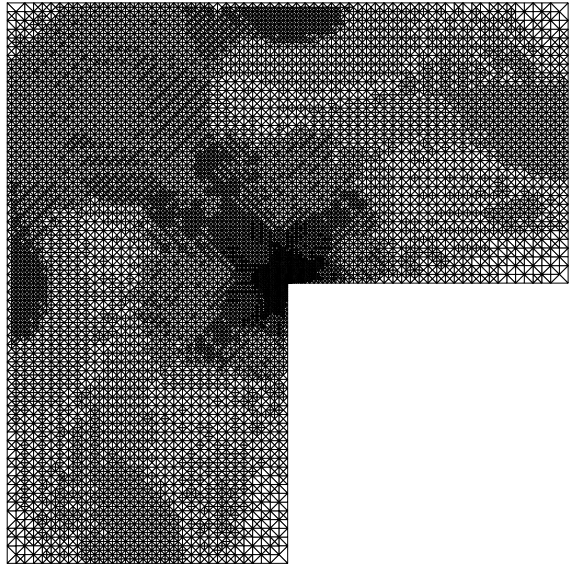}
 \end{minipage}
\caption{Adapted mesh  for the second eigenvalue}
\end{figure}
\begin{figure}
 \begin{minipage}[c]{0.42\textwidth}
\includegraphics[width=4in]{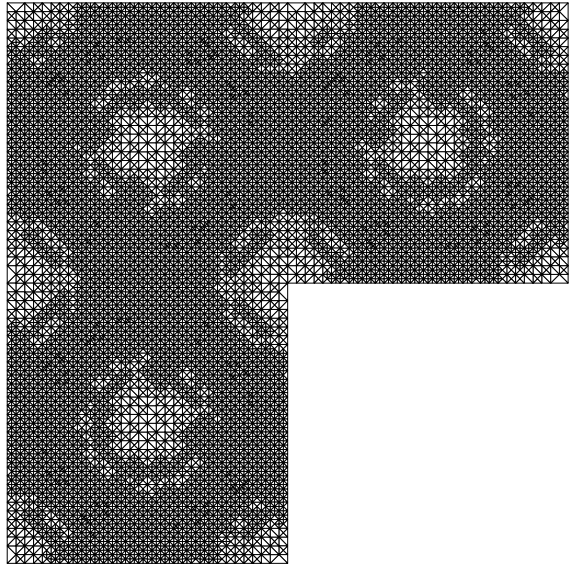}
 \end{minipage}
\caption{Adapted mesh  for the third eigenvalue}
\end{figure}
\newpage
\begin{figure}\begin{center}
\includegraphics[width=3in]{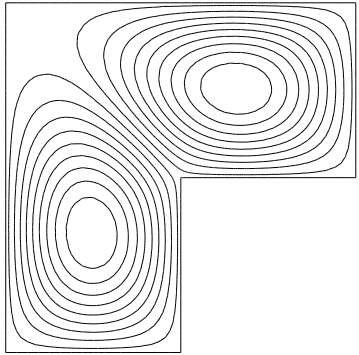}
\includegraphics[width=3in]{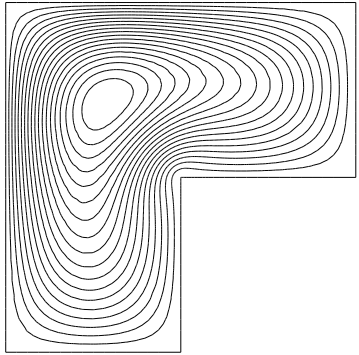}
\includegraphics[width=3in]{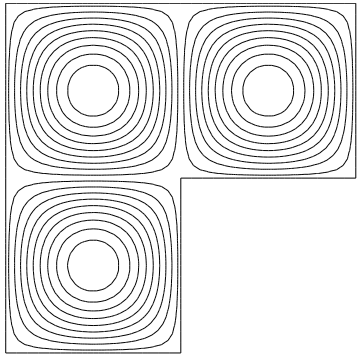}
\end{center}
\caption{The first three eigensolutions}
\end{figure}
\begin{figure}
\begin{center}
\includegraphics[width=3.5in]{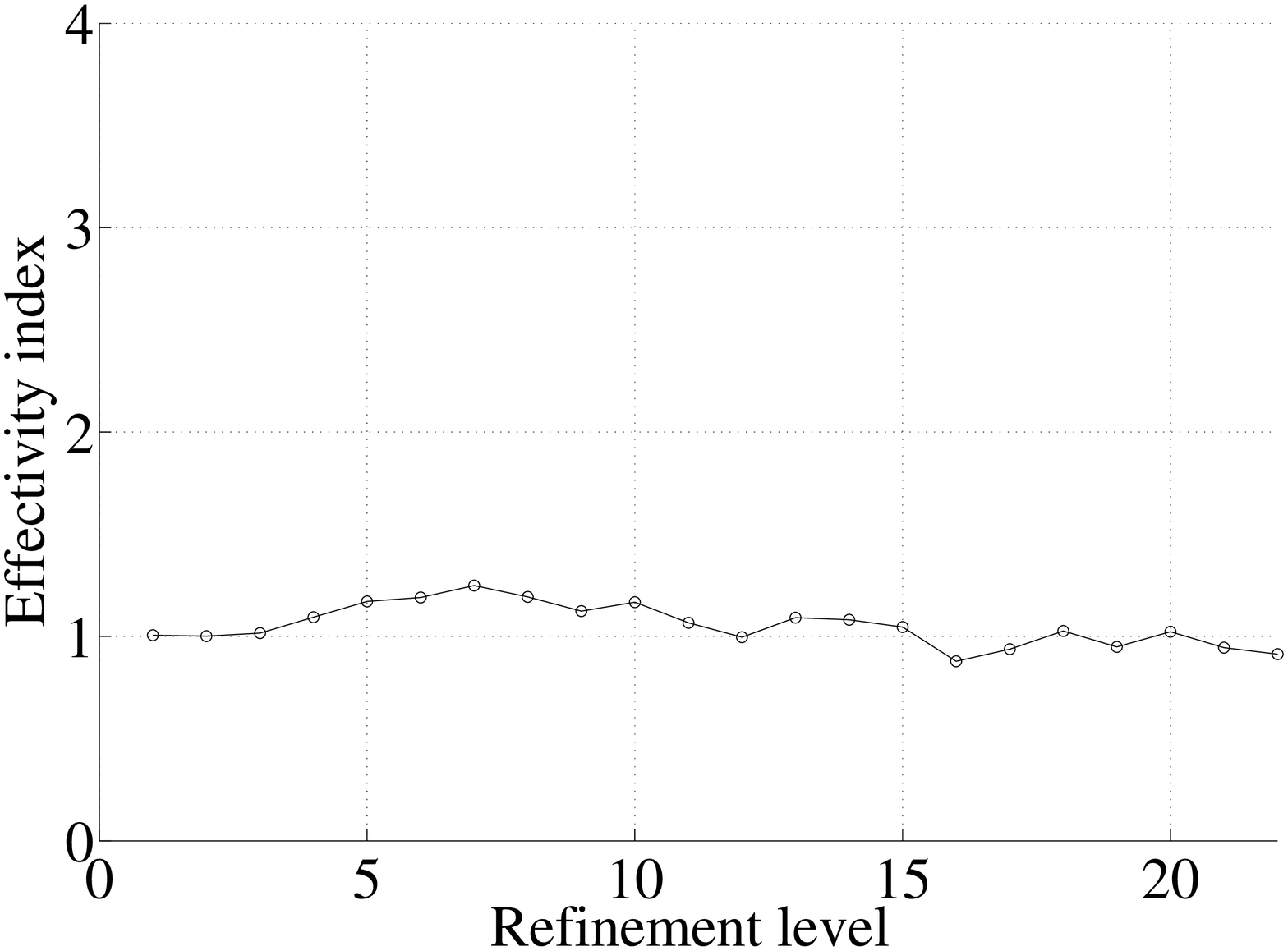}
\includegraphics[width=3.5in]{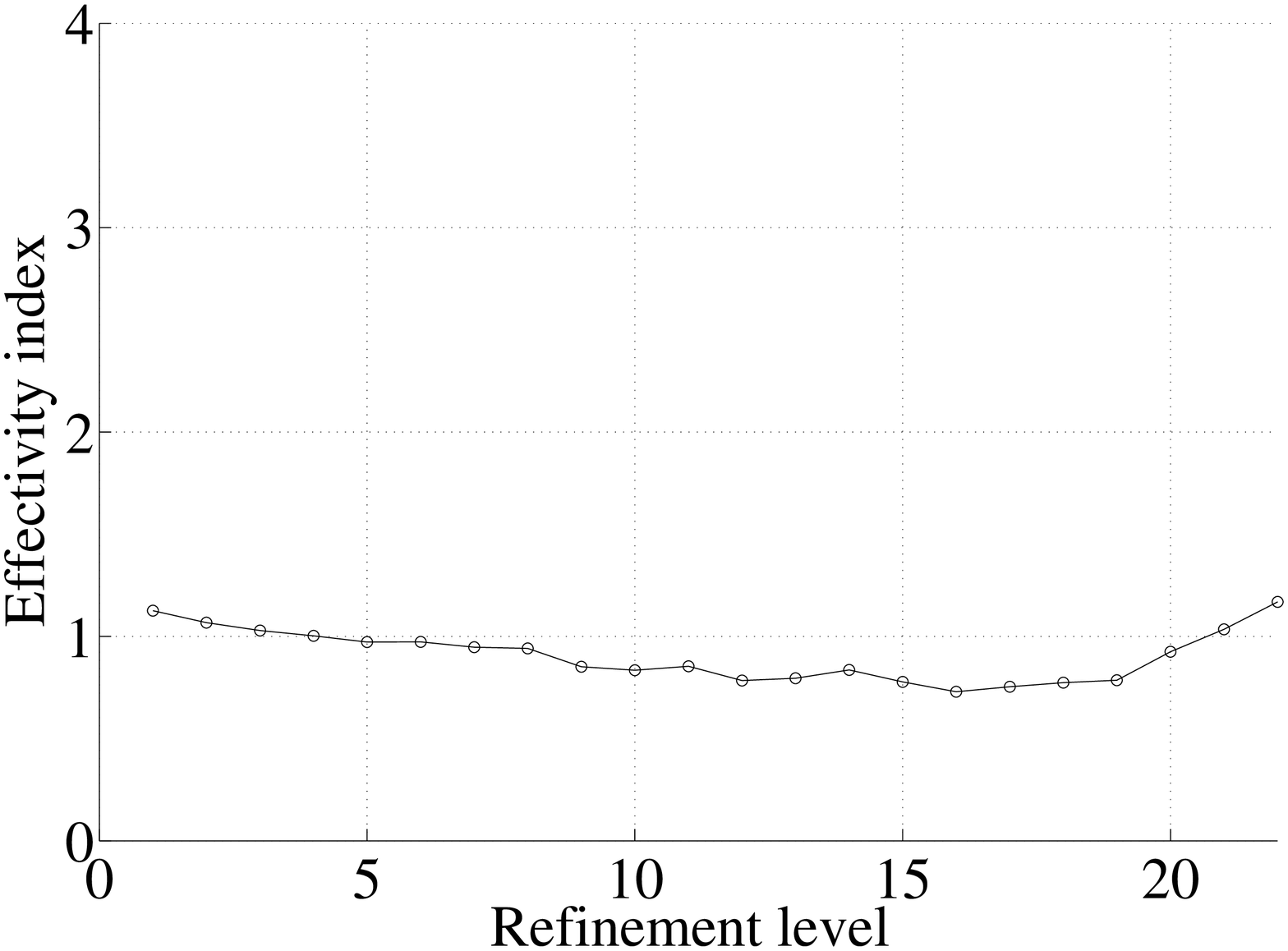}
\includegraphics[width=3.5in]{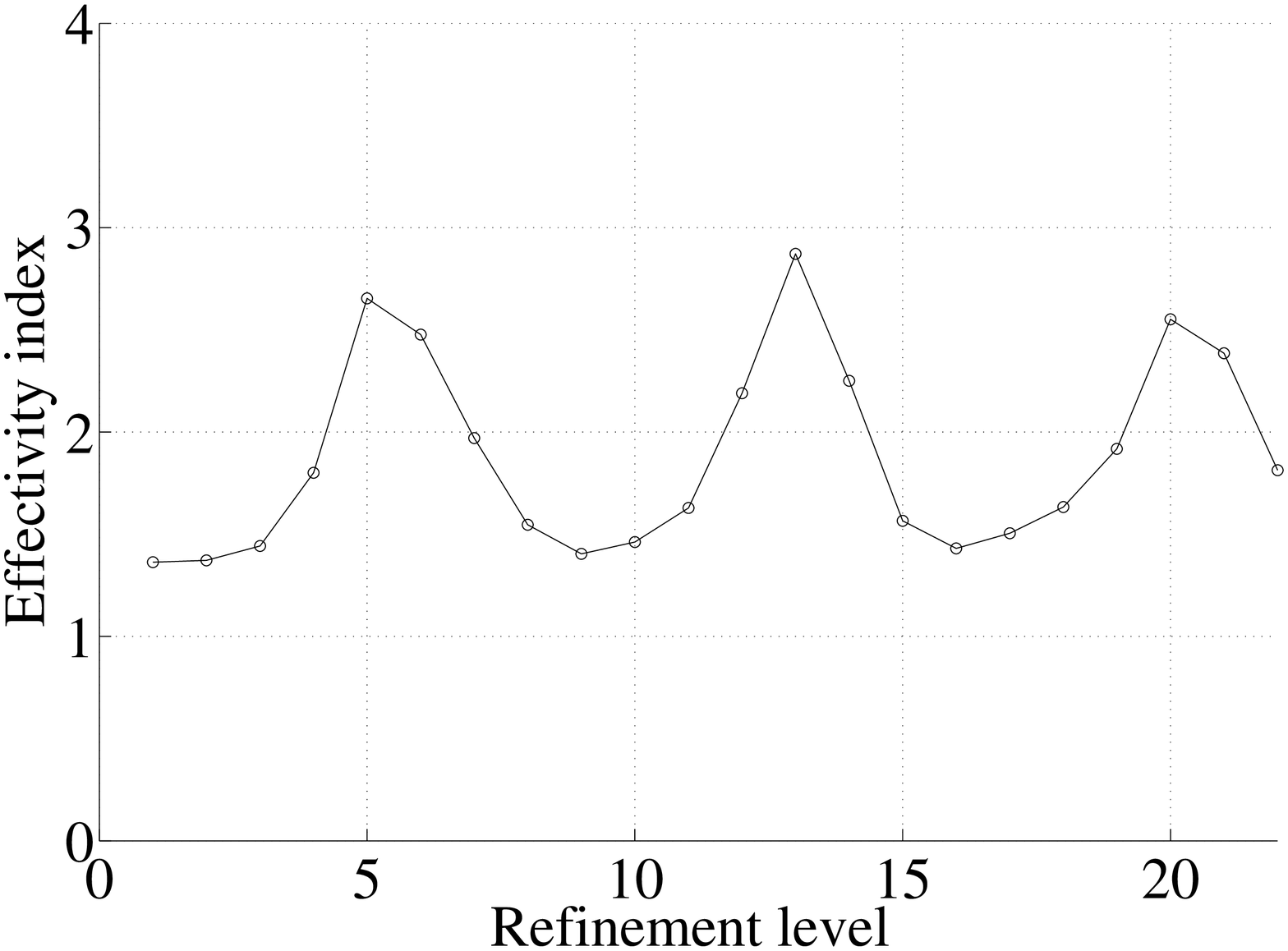}
\end{center}
\caption{Computed effectivity indices for the first three eigenvalue computations.}
\end{figure}
\begin{figure}
\begin{center}
\includegraphics[width=4in]{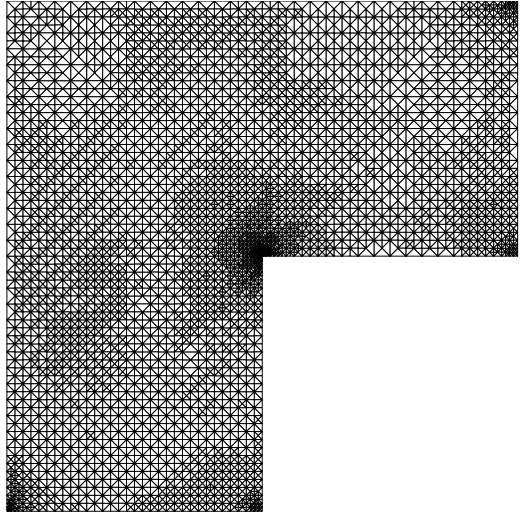}
\end{center}
\caption{Adapted mesh using the full estimate.}
\end{figure}
\begin{figure}
\begin{center}
\includegraphics[width=4in]{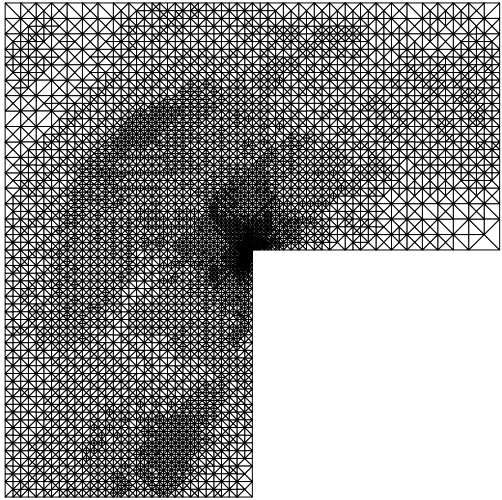}
\end{center}
\caption{Adapted mesh for a partial estimate (only the plate residual).}
\end{figure}
\begin{figure}
\begin{center}
\includegraphics[width=4in]{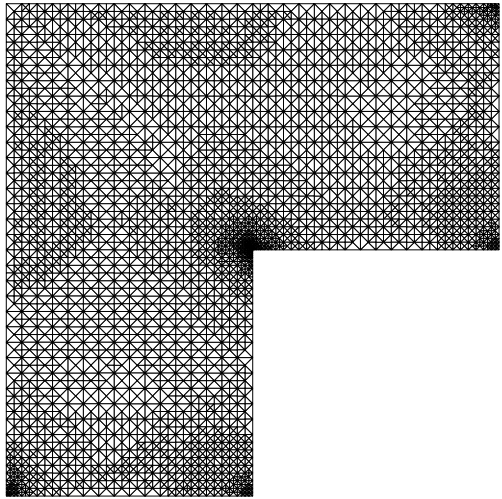}
\end{center}
\caption{Adapted mesh for a partial estimate (only the stress residual).}
\end{figure}
\begin{figure}
\begin{center}
\includegraphics[width=4in]{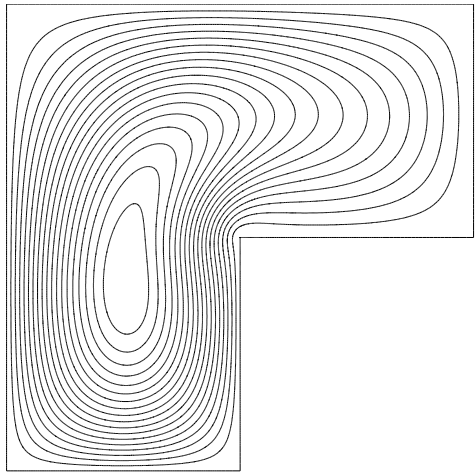}
\end{center}
\caption{Buckling mode.}
\end{figure}
\begin{figure}
\begin{center}
\includegraphics[width=2.5in]{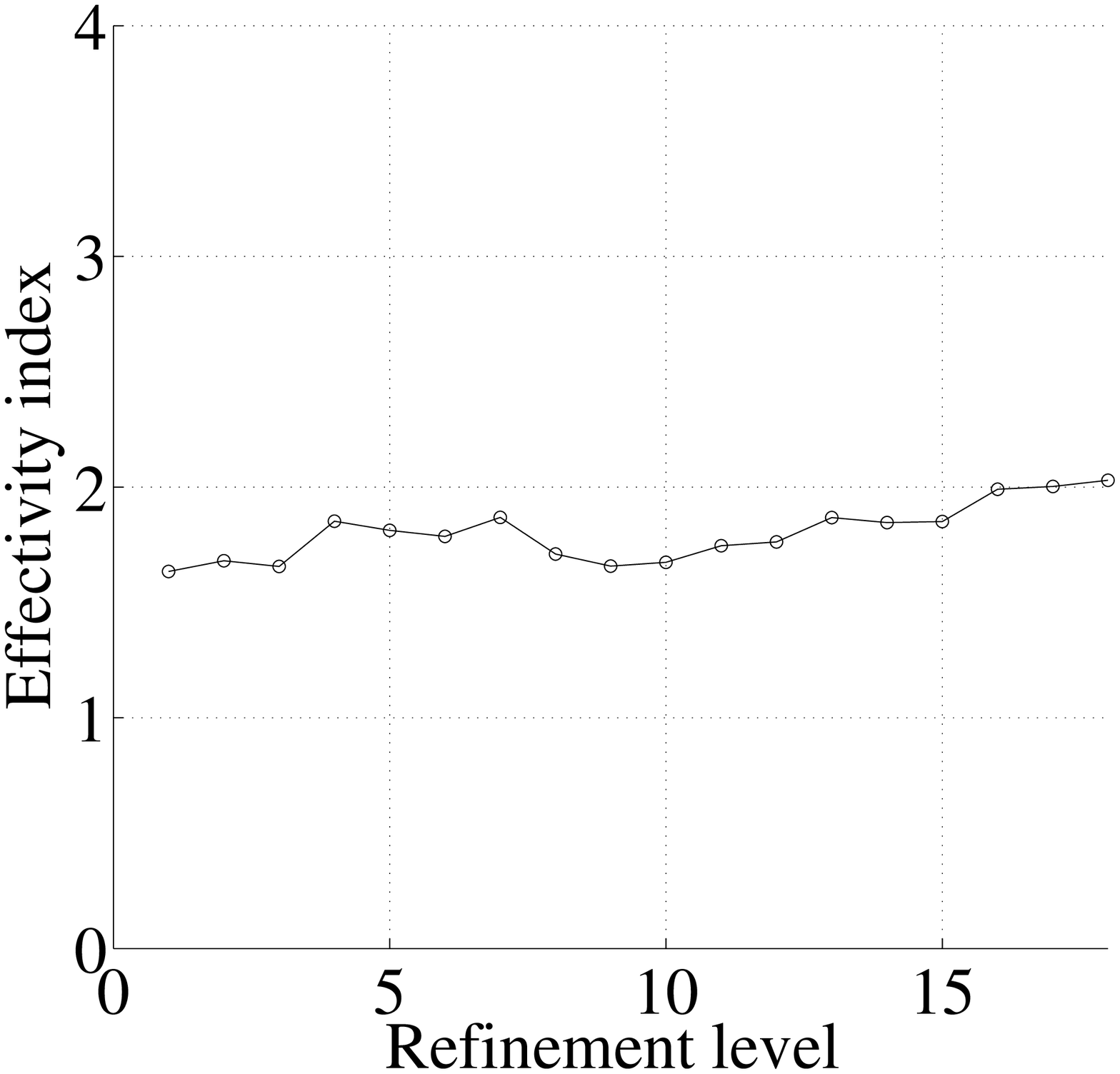}
\includegraphics[width=2.5in]{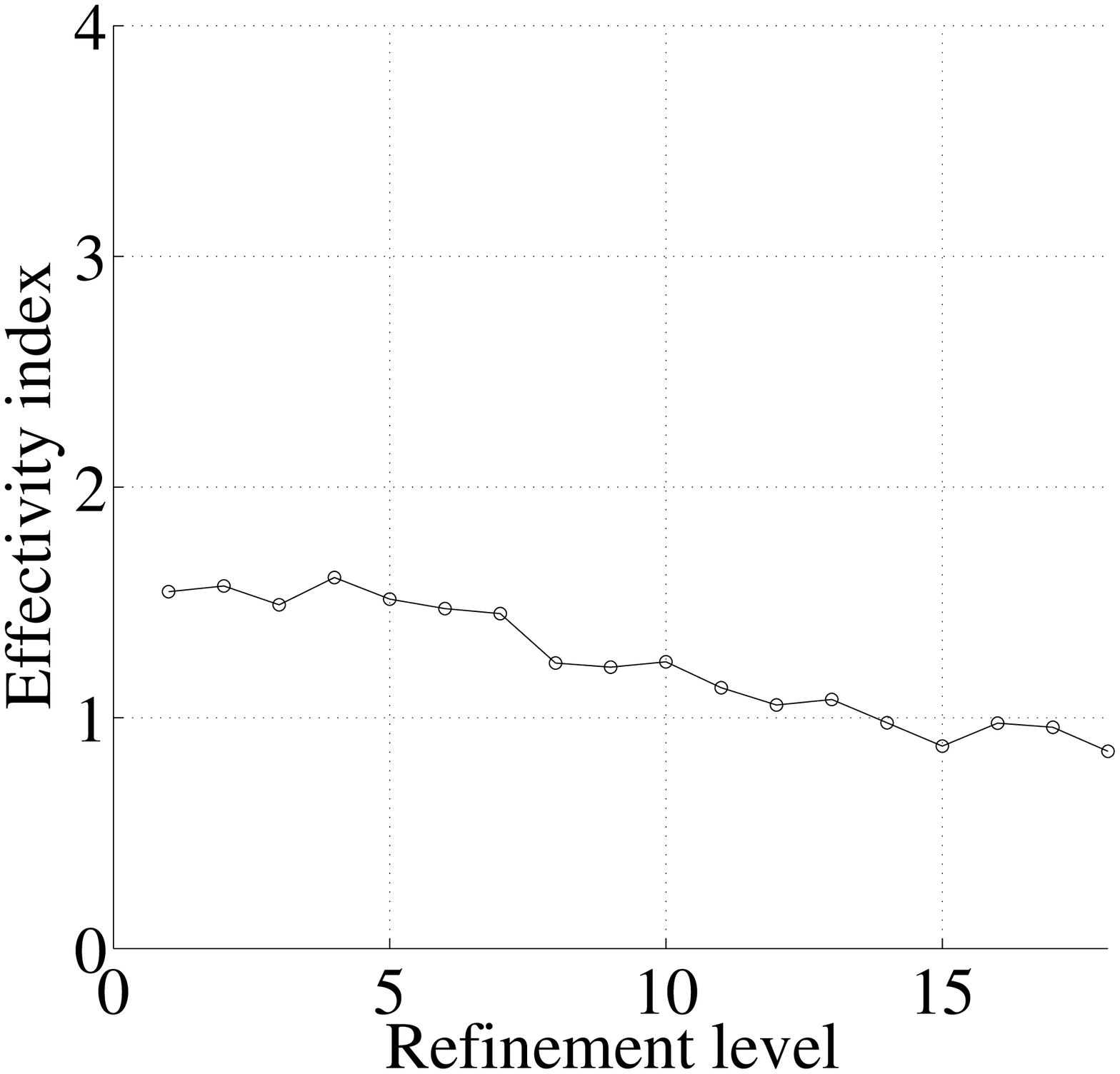}
\includegraphics[width=2.5in]{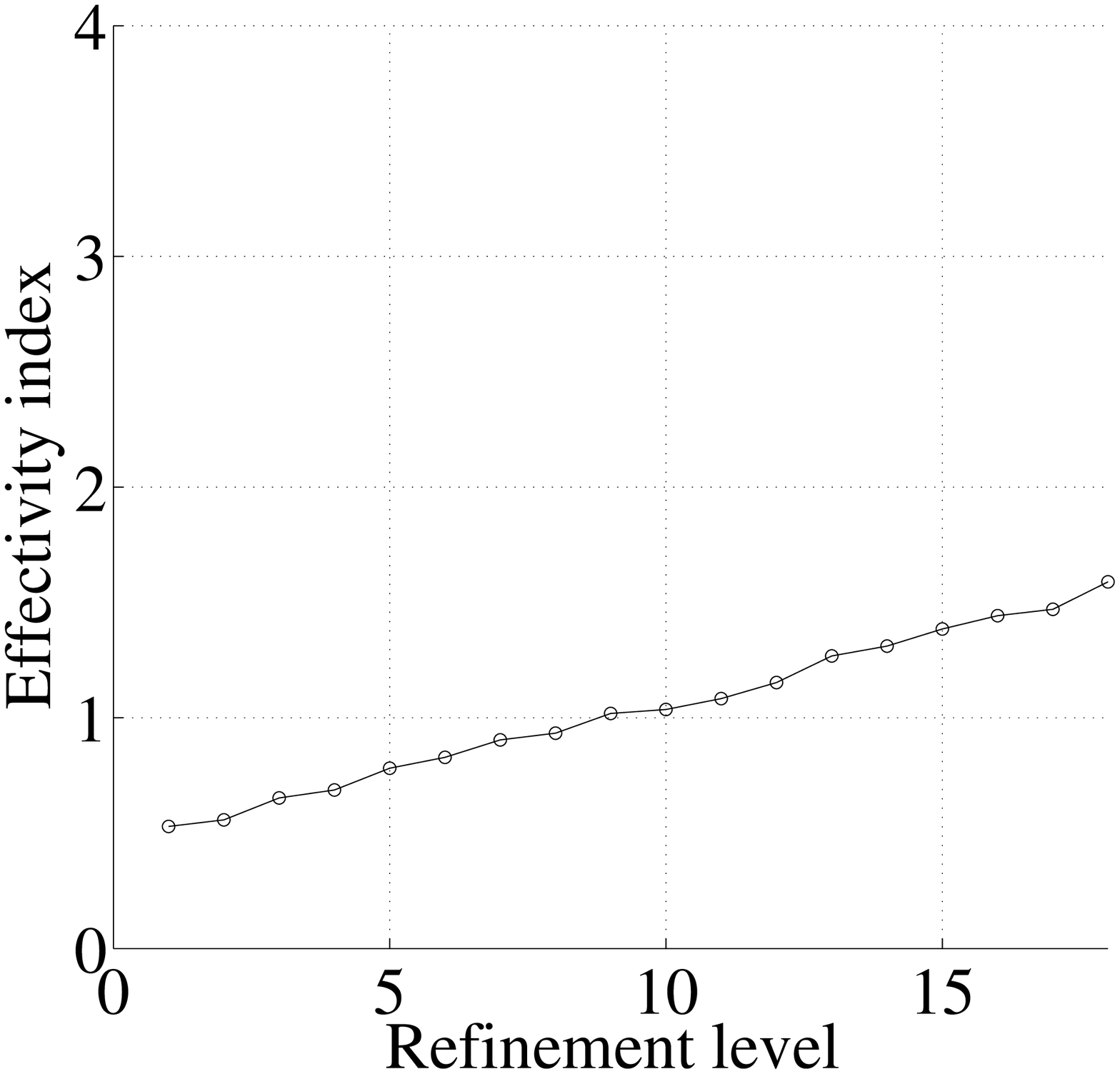}
\end{center}
\caption{Computed effectivity indices for the full estimate and the partial estimates (plate residual and stress residual).}
\end{figure}

\end{document}